\documentclass[11pt]{amsart}

\usepackage[style=alphabetic]{biblatex} 
\usepackage{amsmath}
\usepackage{amsthm}
\usepackage{amsfonts}
\usepackage{amssymb}
\usepackage[all]{xy}
\usepackage{alltt}
\usepackage{enumitem}
\usepackage{stmaryrd}
\usepackage{comment}
\usepackage{xspace}
\usepackage{comment}
\usepackage{libertine}
\usepackage[libertine]{newtxmath}

\theoremstyle{plain}
\newtheorem{proposition}{Proposition}[section]

\newtheorem{theorem}[proposition]{Theorem}
\theoremstyle{definition}

\newtheorem{example}[proposition]{Example}

\newcommand{\Z}{\ensuremath{\mathbf{Z}}}
\newcommand{\F}{\ensuremath{\mathbf{F}}}

\newcommand{\Q}{\ensuremath{\mathbf{Q}}}

\newcommand{\C}{\ensuremath{\mathbf{C}}}

\DeclareMathOperator{\Hom}{Hom}

\DeclareMathOperator{\End}{\mathrm{End}}

\renewcommand{\F}{\ensuremath{\mathbf{F}}}
\renewcommand{\Z}{\ensuremath{\mathbf{Z}}}

\usepackage[margin=1.25in]{geometry}

\usepackage[
        colorlinks=true,urlcolor=blue,linkcolor=blue,citecolor=blue
]{hyperref}

\newcommand{\gbu}[1]{\ensuremath{[#1]}}

\begin{document}  
 
\title{Fuchs' problem for endomorphisms of abelian groups}
\date{\today}
 
\author{Sunil K. Chebolu}
\address{Department of Mathematics \\
Illinois State University \\
Normal, IL 61790, USA}
\email{schebol@ilstu.edu}

\author{Keir Lockridge} 
\address {Department of Mathematics \\
Gettysburg College \\
Gettysburg, PA 17325, USA}
\email{klockrid@gettysburg.edu}

\thanks{Sunil Chebolu is partially supported by Simons Foundation's Collaboration Grant for Mathematicians (516354)}

 
\begin{abstract}
L\'{a}szl\'{o} Fuchs posed the following question: which abelian groups arise as the group of units in a ring? In this paper, we investigate a related question: for such realizable groups $G$, when is there a ring $R$ with unit group $G$ such that every group endomorphism of $G$ is induced by a ring endomorphism of $R$? We answer this question for four common classes of groups: torsion-free abelian groups, groups of odd order, torsion abelian groups, and finitely generated abelian groups.
\end{abstract}
 
\maketitle
\thispagestyle{empty}

\tableofcontents


\section{Introduction}

In the 1960s, L\'{a}szl\'{o} Fuchs posed the problem of classifying the abelian groups that occur as the group of units in some ring (\cite[Problem 72]{F60}). A (possibly nonabelian) group $G$ is called {\bf realizable} if there is a ring $R$ with identity whose unit group $R^\times$ is isomorphic to $G$; in this situation, we say that $R$ {\bf realizes} $G$. Though Fuchs' problem remains open, many authors have made significant progress. For example, the realizable groups have been classified in the following families: cyclic groups (\cite{Gil63, PS70}), groups of odd order (\cite{Ditor71}), alternating and symmetric groups (\cite{DO14b}), finite simple groups (\cite{DO14a}), indecomposable abelian groups (\cite{CL15}), and dihedral groups (\cite{CL17d}). Classifying the realizable 2-groups seems particularly challenging, though partial results appear in \cite{CL17} and \cite{SW20}. Rather than examine the problem from the perspective of the group, other authors have chosen to determine the realizable groups in certain classes of rings; for example, in \cite{D20} (following \cite{DD18a} and \cite{DD18b}) the author determines the finitely generated abelian groups that occur as the group of units in an integral domain, torsion-free ring, or reduced ring. Questions in the spirit of Fuchs' problem touch upon interesting ideas in ring theory, group theory, field theory, and number theory. For instance, they led to new characterizations of Mersenne and Fermat primes; see \cite[p. 328]{CL15}.

In this paper we consider a closely related variant of Fuchs' problem concerning the realizability of group endomorphisms. Let $\mathrm{End}(X)$ denote the set of ring or group endomorphisms of $X$ (depending on whether $X$ is a group or a ring). The units functor $R \mapsto R^\times$ induces a map \[ \xymatrix{ \mathrm{End}(R) \ar[r] & \mathrm{End}(R^\times)}\] (also denoted by $f \mapsto f^\times$) because ring homomorphisms send units to units. We will refer to endomorphisms in the image of this map as {\bf realizable}. If $G$ is realizable and the above map is surjective, then we will say that $G$ is \textbf{fully realizable} and that $R$ \textbf{fully realizes $G$}. We chose this terminology to emphasize the connection to fullness for functors: a functor between locally small categories is full if it is surjective on homsets.  Our goal is to classify the fully realizable groups. In this paper we find all fully realizable groups in the following classes: torsion-free abelian groups, groups of odd order, torsion abelian groups, and finitely generated abelian groups. In subsequent work, the authors will address the full realizability question for several classes of nonabelian groups, including dihedral, (generalized) quaternion, finite simple, symmetric, alternating, and finite general linear groups (\cite{CL}).

 Fuchs' original question may be reframed in light of our definition of fully realizable: a group $G$ is realizable if and only if the identity endomorphism of $G$ is induced by an endomorphism of some ring $R$ via the units functor. A group $G$ might not be realizable at all or---if it {\it is} realizable---the identity endomorphism of $G$ might be the only realizable endomorphism of $G$. For example, $\C_{16}$ (the cyclic group of order 16) is realizable only in characteristics 17 and 34 (see \cite{CL15}); hence, since any endomorphism of a ring $R$ fixes the additive subgroup of $R$ generated by $1$, any endomorphism of a ring realizing $\C_{16}$ will fix $\C_{16}$ element-wise. Our question lies at the other extreme, in that we require that {\em all} endomorphisms be realizable. One could also ask about groups satisfying some intermediate condition. For example, which realizable groups $G$ are realized by a ring that realizes every automorphism of $G$? One can show that if $G$ is realizable, then there is a ring that realizes every inner automorphism of $G$, so any realizable group without outer automorphisms would satisfy this property.

One might expect the fully realizable groups to be quite rare, and this indeed seems to be the case. For example, we will see that the only cyclic groups that are fully realizable are the cyclic groups $\mathbf{C}_n$ for $n \mid 12$ and the infinite cyclic group $\mathbf{C}_\infty$. On the other hand, every torsion-free abelian group $A$ is fully realized by the group algebra $\mathbf{F}_2[A]$ (see the discussion in the first paragraph of \S 3, which establishes item (\ref{tfag}) below). The following theorem summarizes our main results. 

\begin{theorem}
\phantom{Let $G$ be a group.}
\begin{enumerate}
    \item Every torsion-free abelian group is fully realizable. \label{tfag}
    \item The fully realizable groups of odd order are $\C_1$ and $\C_3$.
    \item The fully realizable torsion abelian groups are the groups of the form $W \times H$ where $W$ is an elementary abelian $2$-group and $H$ is a subgroup of $\C_{12}$.
    \item The fully realizable finitely generated abelian groups are the groups $W \times H$ and $W \times \C_\infty^n$ ($n \geq 1$) where $W$ is a finite elementary abelian $2$-group and $H$ is a subgroup of $\C_{12}.$
\end{enumerate}
\end{theorem}

The remainder of this paper is organized as follows. In \S 2 we begin with a few simple but illustrative examples. In \S 3 we prove several general facts, many of which apply more broadly to the category of groups (not just to abelian groups). In the next three sections, we complete the proof of the above theorem, considering groups of odd order in \S 4, torsion abelian groups in \S 5, and finitely generated abelian groups in \S 6.

\vspace{.1in}

\noindent {\bf Acknowledgement.} We wish to thank the anonymous referee for carefully reading our paper, making detailed comments, and encouraging us to improve the proof of Proposition \ref{c3iq}.

\section{Basic Examples}

In this section, we occasionally make reference to the fact that if a group is fully realizable, then it is fully realized by a ring of characteristic 2 that is generated by its units, and it is not fully realizable in any other characteristic; see Proposition \ref{basics}. 

\begin{example}[The cyclic group of order 2 is fully realizable]
The group $\mathbf{C}_2$ is fully realized by the group algebra $\mathbf{F}_2[\mathbf{C}_2]$; the identity map induces the identity homomorphism and the augmentation map followed by the inclusion of $\F_2$ into $\F_2[\C_2]$ induces the trivial homomorphism.
\end{example}

\begin{example}[The cyclic group of order 4 is fully realizable]
Write $\mathbf{C}_4 = \langle x \rangle.$ The quotient $\mathbf{F}_2[\mathbf{C}_4]/((x+1)^3) \cong \mathbf{F}_2[z]/(z^3)$ has unit group $\mathbf{C}_4$. There are 4 ring endomorphisms determined by the assignments $z \mapsto z, z \mapsto z^2, z \mapsto z + z^2,$ and $z \mapsto 0$. These maps induce the four group endomorphisms of $\mathbf{C}_4$: $x \mapsto x$, $x \mapsto x^2$, $x \mapsto x^3$ and $x \mapsto 1$.
\end{example}

\begin{example}[The groups $\mathbf{C}_{2^k}$ for $k \geq 3$ are not fully realizable]
For $k \geq 3$, the group $\mathbf{C}_{2^k}$ is not realizable in characteristic 2 (\cite{CL15}) and hence it is not fully realizable.
\end{example}

\begin{example}[The cyclic group of order 3 is fully realizable]
The group $\mathbf{C}_3 = \langle x \rangle$ is fully realizable, but it is not fully realized by every ring $R$, generated by its units, with $R^\times = \mathbf{C}_3$. For example, $\mathbf{F}_4^\times = \mathbf{C}_3$, but every ring endomorphism of $\mathbf{F}_4$ is also a field automorphism; hence, only the maps $x \mapsto x$ and $x \mapsto x^2$ are realizable. The trivial homomorphism is not induced by a ring endomorphism of $\mathbf{F}_4$. Now consider the ring $R = \F_2[\C_3] = \mathbf{F}_2 \times \mathbf{F}_4$; it is also generated by its units and has unit group $\mathbf{C}_3$. The map $(a, b) \mapsto (a, b)$ realizes the identity map, the map $(a, b) \mapsto (a, b^2)$ realizes $x \mapsto x^2$, and the map $(a, b) \mapsto (a, a)$ realizes the trivial homomorphism. Thus, $R$ fully realizes $\mathbf{C}_3$. We will later prove that no other nontrivial group of odd order is fully realizable.
\end{example}

\begin{example}[Torsion-free abelian groups are fully realizable]
Concretely, the infinite cyclic group $\mathbf{C}_\infty = \langle x \rangle$ is fully realized by $R = \mathbf{F}_2[x^{\pm1}]$ because every ring endomorphism of $R$ determines, and is determined by, the image of $x$, which can map to any unit. Every map from $\C_\infty$ to a realizable group is realized by a ring homomorphism. As mentioned in the introduction, for any torsion-free abelian group $A$, the group algebra $\mathbf{F}_2[A]$ has group of units $A$, and hence it fully realizes $A$ (see \S 3 for more details).
\end{example}

Kaplansky's unit conjecture entails that $\mathbf{F}_2[A]^\times = A$ for {\em all} torsion-free groups $A$; this conjecture was recently proved false by Gardam in \cite{Gardam21} (though it still may be the case that some nonabelian torsion-free group is fully realizable).

\begin{example}[The cyclic groups of orders 6 and 12 are fully realizable] \label{psrg} A product of fully realizable groups may or may not be fully realizable. For the group $\mathbf{C}_2 \times \mathbf{C}_3$ the obvious trick works:  the product of  rings $\mathbf{F}_2[\mathbf{C_2}] \times (\mathbf{F}_2 \times \mathbf{F}_4)$, where the factors fully realize their unit groups, fully realizes $\mathbf{C}_2 \times \mathbf{C}_3$. This is because, since there are no maps from $\mathbf{C}_2$ to $\mathbf{C}_3$ or vice-versa, every endomorphism of $\mathbf{C}_2 \times \mathbf{C}_3$ has the form $(\phi_1, \phi_2)$ where $\phi_1 \in \mathrm{End}(\mathbf{C}_2)$ and $\phi_2 \in \mathrm{End}(\mathbf{C}_3).$ We may reason similarly to conclude that $\mathbf{C}_4 \times \mathbf{C}_3$ is fully realizable.
\end{example}

\begin{example}[The group $\mathbf{C}_3 \times \mathbf{C}_3$ is not fully realizable] \label{c3xc3}
If this group is the group of units in a ring $R$ of characteristic 2, generated by its units, then $R$ is a finite product of fields (\cite{Ditor71}). Given that $R$ is generated by its units, we must have that $R$ is $\mathbf{F}_2^a \times \mathbf{F}_4^2,$ where $a \in \{ 0, 1\}.$ As we will later show, this ring has $16$ or $25$ endomorphisms (depending on the value of $a$), but $\mathbf{C}_3 \times \mathbf{C}_3$ has 81 endomorphisms, so not every group endomorphism is realizable.
\end{example}

We will generalize the situation in Example \ref{psrg} concerning when the product of fully realizable groups is actually fully realizable; see Proposition \ref{products}. We will also later prove that if a group is fully realizable, then so are its direct summands; see Proposition \ref{summand}.

\section{Preliminaries}

Given a ring $R$ and a subgroup $H \leq R^\times$, we will write $\gbu{H}$ for the subring of $R$ generated by $H$. A ring $R$ is generated by its units if $R = \gbu{R^\times}.$ For any field $k$ and group $G$, let $k[G]$ denote the group algebra of $G$ over $k$. Any group homomorphism $\phi \colon G \longrightarrow H$ determines a unique ring homomorphism $\overline{\phi} \colon k[G] \longrightarrow k[H]$ such that $\overline{\phi}|_G = \phi.$ This assignment induces an injective map \[\mathrm{Hom}_{\mathrm{Groups}}(G, H) \longrightarrow \mathrm{Hom}_{\mathrm{Rings}}(k[G], k[H])\] that is also surjective if $(k[H])^\times = H$. In particular, any group $G$ with $(\mathbf{F}_2[G])^\times = G$ is fully realizable (e.g., if $G$ is a torsion-free abelian group; see \cite[Theorem 3.3]{CL15}). 

When $R$ fully realizes $G$ and $\phi \in \End(G)$, we will let $\hat{\phi}$ denote an arbitrarily chosen element of $\End(R)$ with $\hat{\phi}^\times = \phi.$ The next proposition implies that in order to classify the fully realizable groups, it suffices to consider rings of characteristic 2 that are generated by their units. Such a ring $R$ has the form $\mathbf{F}_2[R^\times]/I$ for some ideal $I \subseteq \mathbf{F}_2[R^\times]$ (up to isomorphism). 

\begin{proposition} \label{basics}
Let $R$ be a ring and let $G$ be a group.
\begin{enumerate}
\item If $R$ fully realizes $G$, then $R$ has characteristic 2 and the subring of $R$ generated by its units fully realizes $G$. \label{gen-by-units}
\item If $R$ has characteristic $2$, is generated by its units, and realizes $G$, then $R \cong \mathbf{F}_2[G]/I$ for some (two-sided) ideal $I$ of the group algebra $\mathbf{F}_2[G].$ \label{gaq}
\end{enumerate}
\vspace{0.05in}
\end{proposition}
\begin{proof}
If $R$ is a ring whose characteristic is not 2, then the unit $-1$ is not the trivial unit $1$, and any ring endomorphism of $R$ must fix $-1$. Hence, the trivial homomorphism is not realizable and $R$ does not fully realize its unit group. Thus, if a group $G$ is fully realizable, then it must be fully realizable only in characteristic 2. Now suppose that $R$ fully realizes $R^\times$ (hence, $\mathrm{char}\, R = 2$). Let $S$ be the subring of $R$ generated by its units; note that $S^\times = R^\times$. Take $\phi \in \mathrm{End}(R^\times).$ Since any ring endomorphism of $R$ must take units to units, the restriction of $\hat{\phi}\in \End(R)$ to $S$ gives a map $\psi \in \End(S)$ with $\psi^\times = \hat{\phi}^\times = \phi$. Thus, $S$ fully realizes $S^\times$. 

If $R$ has characteristic 2, is generated by its units, and realizes $G$, then the inclusion of $G$ into $R$ induces a surjective ring homomorphism $\mathbf{F}_2[G] \longrightarrow R.$ Item (\ref{gaq}) follows.
\end{proof}

The `if' direction of the following proposition provides a strategy for proving that a group is fully realizable.

\begin{proposition} \label{ideal-pres}
Suppose $R$ is a ring of characteristic 2 that is generated by its group of units $G$. Let $I$ denote the kernel of the surjective ring homomorphism \[ \pi\colon \F_2[G] \longrightarrow R\] induced by the inclusion of $G$ into $R$. The ring $R \cong \F_2[G]/I$ fully realizes $G$ if and only if $\overline{\phi}(I) \subseteq I$ for all $\phi \in \mathrm{End}(G).$
\end{proposition}
\begin{proof}
Suppose $R$ fully realizes $G$ and take $\phi \in \mathrm{End}(G)$. Since $\pi \circ \overline{\phi} = \hat{\phi} \circ \pi$ on $G$, the following diagram commutes: \[ \xymatrix{\F_2[G] \ar[r]^-{\overline{\phi}} \ar[d]_{\pi}& \F_2[G] \ar[d]^-\pi \\ R \ar[r]_-{\hat{\phi}} & R. }\] This implies $\overline{\phi}(I) \subseteq I.$

Conversely, suppose $\overline{\phi}(I) \subseteq I$ for all $\phi \in \mathrm{End}(G)$ and fix $\phi \in \End(G)$. Since $\pi(\overline{\phi}(I)) \subseteq \pi(I) = 0$, there is a map $\psi \in \End(R)$ such that $\psi \circ \pi = \pi \circ \overline{\phi}$, and \[\psi^\times(g) = \psi(\pi(g)) = \pi(\overline{\phi}(g)) = \pi(\phi(g)) = \phi(g)\] (recall that $\pi$ includes $G$ into $R$). Thus, $R$ fully realizes $G$.
\end{proof}

The next proposition implies that summands of fully realizable groups are fully realizable.

\begin{proposition} \label{summand}
Suppose $G$ is a group with an idempotent endomorphism $\phi$ and let $H$ denote the image of $\phi$. If $R$ is a ring, generated by its units, that fully realizes $G$, then the subring $[H] \subseteq R$ fully realizes $H$. In particular, any direct summand of a fully realizable group is fully realizable.
\end{proposition}
\begin{proof}
Define $G$, $\phi$, and $H$ as in the statement of the theorem and suppose $R$ is a ring, generated by its units, that fully realizes $G$. The map $\phi$ lifts to an idempotent endomorphism $\hat{\phi}$ of $R$. Let $S = [H]$. Since $\hat{\phi}$ fixes $S$ element-wise, $S$ cannot contain an element of $G \setminus H$, thus $S^\times = H$. If $\rho \in \mathrm{End}(H),$ then the map $\rho' \in \mathrm{End}(G)$ defined by $\rho'(g) = \rho(\phi(g))$ is induced by a map $K \in \mathrm{End}(R)$ (that is, $K^\times = \rho'$). Further, $K|_S$ maps $S$ to $S$ and $(K|_S)^\times = \rho'|_S = \rho.$ Thus, $S$ fully realizes $H$.
\end{proof}

We next describe a situation where the product of fully realizable groups is again fully realizable. It applies, for example, when the pair of groups have relatively prime orders.

\begin{proposition} \label{products}
If $G$ and $H$ are fully realizable groups such that there are no nontrivial homomorphisms from $G$ to $H$ or from $H$ to $G$, then $G \times H$ is fully realizable.
\end{proposition}
\begin{proof} For any groups $G, H$ and $C$, $\mathrm{Hom}(G \times H, C)$ is in one-to-one correspondence with pairs \[ (\tau, \sigma) \in \mathrm{Hom}(G, C) \times \mathrm{Hom}(H, C) \] such that \[ \tau(x)\sigma(y) = \sigma(y)\tau(x)\] for all $x \in G$ and $y \in H$. Hence, if there are no nontrivial homomorphisms from $G$ to $H$ or from $H$ to $G$, then we have \[\mathrm{End}(G \times H) \cong \mathrm{End}(G) \times \mathrm{End}(H).\]

Now suppose $R$ and $S$ fully realize $G$ and $H$ (as in the statement of the proposition). By hypothesis, any $\phi \in \mathrm{End}(G \times H)$ has the form $\phi = \phi_1\times \phi_2$ for $\phi_1 \in \mathrm{End}(G)$ and $\phi_2 \in \mathrm{End}(H).$ There are endomorphisms $F_1 \in \mathrm{End}(R)$ and $F_2 \in \mathrm{End}(S)$ with $F_1^\times = \phi_1$ and $F_2^\times = \phi_2$. We now have $F_1 \times F_2 \in \mathrm{End}(G \times H)$ with $(F_1 \times F_2)^\times = \phi_1\times \phi_2 = \phi.$ Thus, $G \times H$ is fully realizable.\end{proof}

\begin{example} \label{2torsionxc3}
For any fully realizable $2$-group $T$, $T \times \mathbf{C}_3$ is fully realizable.
\end{example}

Next, we give a useful tool that says certain products of rings cannot fully realize their groups of units. If $G$ is a group, let $Z(G)$ denote the center of $G$.

\begin{proposition} \label{products-fail}
Let $R$ and $S$ be rings of characteristic 2. Suppose that, for any decomposition $S = A \times B$, if $A \neq \{0\}$ then $A^\times$ is nontrivial. If there exists a nontrivial group homomorphism $\phi\colon R^\times \longrightarrow Z(S^\times)$, then $R \times S$ does not fully realize its group of units.
\end{proposition}
\noindent Note that the hypothesis on $S$ is satisfied if $S$ is indecomposable as a ring with $S^\times \neq \{1\}$ or if $S$ is an algebra over $\F_{2^k}$ for $k > 1$.
\begin{proof} The map $\psi(u, v) = (u, \phi(u)v)$ defines an endomorphism of $R^\times \times S^\times$ since the image of $\phi$ is a subset of $Z(S^\times).$ We claim $\psi$ cannot extend to a ring endomorphism of $R \times S$. Assume to the contrary that such an extension, $\overline{\psi}$, exists. Pick any $u \in R^\times$ such that $\phi(u) \neq 1$. Now, $\psi(u, 1) = (u, \phi(u))$, so \[
\begin{aligned}
(1 + u, 1 + \phi(u)) &= \overline{\psi}(1 + u, 0)\\
&= \overline{\psi}((1, 0)\cdot (1 + u, 0)) \\
&= \overline{\psi}(1,0)\overline{\psi}(1 + u, 0)\\
&= \overline{\psi}(1,0)\cdot(1 + u, 1 + \phi(u)). 
\end{aligned} \] Thus $\overline{\psi}(1,0) = (e, f),$ where $e$ and $f$ are both nonzero central idempotents, and $\overline{\psi}(0,1) = (1 + e, 1 + f).$ Since $f \neq 0$, we have a decomposition $S \cong Sf \times S(1 + f)$ where $Sf \neq \{0\}$. For all $v \in S^\times$, we have
\[
\begin{aligned}
(0, 1+ v) &= (1,1) + (1, v) \\
&= \overline{\psi}(1,1) + \overline{\psi}(1,v) \\
&= \overline{\psi}(0, 1+ v) \\
&= (1 + e, 1 + f)(0, 1 + v).
\end{aligned}
\]
Hence, \[1 + v = (1 + f)(1+v) = 1 + v + f + fv,\] which implies $fv = f$. Thus $Sf$ is a nonzero direct factor of $S$ with a trivial unit group, contradicting our hypothesis on $S$.\end{proof}

\begin{example} \label{cic3ex}
Let $V$ be a torsion-free abelian group with $\C_3$ as a quotient (e.g., $V = \C_\infty$). Consider the ring $R \times S = \mathbf{F}_2[V] \times \mathbf{F}_4$. Since $\mathbf{F}_4$ is indecomposable as a ring and there is a nontrivial map from $V$ to $\mathbf{C}_3$, this particular product ring does not fully realize its unit group $V \times \mathbf{C}_3$. We will later prove that $\mathbf{C}_\infty \times \mathbf{C}_3$ is not fully realizable.
\end{example}

For later use, we conclude this section with relevant background material on group algebras over a field $k$. Let $G$ be an group. The augmentation map of the group algebra $k[G]$ is the ring homomorphism $\epsilon\colon k[G] \longrightarrow k$ that sends each $g \in G$ to 1. The kernel of this map is called the augmentation ideal, and it is generated by elements of the form $1 - g$ for $g \in G$. If $g \in G$ and $G$ appears in a direct product of groups, we will abuse notation and write $g$ for the image of $g$ under the canonical inclusion of $G$ into the direct product. For example, if we write $\C_3 \times \C_4 = \langle x \rangle \times \langle y \rangle$, we write $x$ for $(x, 1)$ and $y$ for $(1, y)$.

\begin{proposition} \label{kgproduct}
Let $k$ be a field with prime subfield $k_0$, let $G_1, \dots, G_n$ be groups, and let \[S = [G_1 \times \cdots \times G_n] \subseteq k[G_1]\times \cdots \times k[G_n].\] Then, \[ S \cong k_0[G_1 \times \dots \times G_n]/I, \]  where $I$ is generated by the elements \[1 - a - b + ab = (1 - a)(1-b)\] for all $a \in G_i$ and $b \in G_j$ when $i \neq j$. In fact, $I$ contains all elements of the form \[\prod_{i=1}^n g_{i} - \sum_{i=1}^ng_{i} + n - 1,\] for all $g_i \in G_i$. Further, if $k = \F_2$ (and hence $k_0 = \F_2$), \[ S^\times = (\F_2[G_1])^\times \times \cdots \times (\F_2[G_n])^\times.\]
\end{proposition}
\begin{proof}
We prove the proposition for $n = 2$ and leave the generalization to the reader. Define the ideal $I$ as in the statement of the theorem. The inclusion of $G \times H$ into $k[G] \times k[H]$ induces a ring homomorphism \[ k_0[G\times H] \longrightarrow k[G] \times k[H]\] whose image is $S$. Take $z = \sum \alpha_i g_i h_i$ in the kernel of this homomorphism; since the image of $z$ is \[(\sum \alpha_i g_i, \sum \alpha_i h_i) = (0, 0),\] each coordinate is zero and has augmentation zero, so $\sum \alpha_i = 0$. Now, \[ \sum \alpha_i(1 - g_i)(1 - h_i) = \sum \alpha_i - \sum \alpha_i g_i - \sum \alpha_i h_i + \sum \alpha_i g_i h_i\] and thus \[z = \sum \alpha_i(1 - g_i)(1 - h_i)= \sum \alpha_i(1 - g_i -h_i + g_i h_i) \in I.\] Conversely, it is easy to check that $I$ is contained in the kernel. This proves that $S = [G\times H] \subseteq k[G] \times k[H]$ is isomorphic to $k_0[G \times H]/I.$

To prove that $S$ has the claimed unit group when $k = \F_2$, it suffices to prove that $(u, 1) \in S^\times$ for all $u \in (\F_2[G])^\times$ (the proof for units of the form $(1, v)$ is similar). We have $u = \sum \alpha_i g_i$ and $\epsilon(u) = \sum \alpha_i = 1$. Now, \[(u, 1) = \sum \alpha_i(g_i, 1) \in S. \] This completes the proof (for $n = 2$).
\end{proof}

\section{Groups of odd order}

In \cite{Ditor71} Ditor proves that if $T$ is a ring whose group of units has odd order and $R$ is the subring of $T$ generated by its units, then $R$ is a finite direct product of fields of characteristic 2 (see also the proof of Proposition 3 in \cite{CL17}). Further, $R$ has at most one factor isomorphic to $\F_2$, for otherwise $R$ would not be generated by its units (for any ring $Q$ and any integer $n \geq 1$, every element $(x, y)$ in the subring of $\F_2^n \times Q$ generated by its units must have $x = 0$ or $x = 1$). Hence, to classify the fully realizable groups of odd order, it suffices to consider rings of the form \[ R = \prod_{i=1}^t \F_{2^{k_i}}, \] where $k_1 \leq k_2 \leq \cdots \leq k_t$ and $k_2 > 1$ if $t > 1$.

\begin{proposition}
The only fully realizable groups of odd order are $\C_1$ and $\mathbf{C}_3$.
\end{proposition}

\begin{proof}
As we have already seen, $\F_2$ fully realizes $\C_1$ and $\F_2 \times \F_4$ fully realizes $\C_3$.

Now suppose $G$ is a fully realizable group of odd order. As explained before the statement of the proposition, we may assume $G$ is fully realized by \[ R = \prod_{i=1}^t \F_{2^{k_i}}, \] where $k_1 \leq k_2 \leq \cdots \leq k_t$ and $k_2 > 1$ if $t > 1$.

If $t = 1$, then $R = \F_{2^k}$ and the group $\mathbf{F}_{2^k}^\times = \mathbf{C}_{2^k -1}$ has $2^k - 1$ endomorphisms. Every ring endomorphism of $\mathbf{F}_{2^k}$ must be a field automorphism, and so there are exactly $k$ ring endomorphisms. Thus, if $\mathbf{F}_{2^k}$ is fully realizable, then $2^k - 1 \leq k$. This forces $k = 1$.

Now suppose $t > 1$ and fix $k = k_i$ for some $k_i > 1$. The subring of $R$ generated by $\F_{2^{k}}^\times$ is isomorphic to $S = \F_2 \times \F_{2^{k}}$. Since $R$ fully realizes its units, by Proposition \ref{summand} we have that $S$ fully realizes its units. However, \[ \End(S) = \Hom_{\mathrm{Rings}}(S, \F_2) \times \Hom_{\mathrm{Rings}}(S, \F_{2^{k}}), \] so $|\End(S)| = 1 \cdot (1 + k) = 1 + k$, and $|\End(\F_{2^k}^\times)| = 2^k - 1$. So we must have \[ k+ 1 \geq 2^k - 1, \] forcing $k = 2.$ Thus, \[ R = \F_2^a \times \F_4^b\] where $a \in \{0, 1\}$ and $b \geq 1$. Since $R^\times = \C_3^b$, in light of Proposition \ref{summand} it now suffices to prove that $G = \C_3 \times \C_3$ is not fully realizable; i.e., we may assume $b = 2$.

We claim there are not enough ring endomorphisms of $R$ to realize all the group endomorphisms of $G$. There are $81$ endomorphisms of $G$. Note that \[ \mathrm{End}(R) = \mathrm{Hom}_{\mathrm{Rings}}(R, \mathbf{F}_2)^a \times \mathrm{Hom}_{\mathrm{Rings}}(R, \mathbf{F}_4)^2.\] We have \[|\mathrm{Hom}_{\mathrm{Rings}}(R, \mathbf{F}_4)| = 4+a\] and \[|\mathrm{Hom}_{\mathrm{Rings}}(R, \mathbf{F}_2)| = a,\] so \[|\mathrm{End}(R)| = (4+a)^2.\] But $81 > 16$ ($a = 0$) and $81> 25$ ($a = 1$). This completes the proof.
\end{proof}

\section{Torsion abelian groups}

In this section, we will prove the following theorem.

\begin{theorem} \label{torsion}
The fully realizable torsion abelian groups are the groups of the form $W \times H$ where $W$ is an elementary abelian $2$-group and $H$ is a subgroup of $\C_{12}$.
\end{theorem}

First, we consider the indecomposable cyclic groups and then prove that a fully realizable torsion abelian group is a direct sum of (fully realizable) indecomposable cyclic groups. Let $\C_{p^\infty}$ denote the $p$-quasicyclic group.

\begin{proposition} Let $p$ be a prime and $n \geq 0$ be an integer or $n = \infty$. The group $\mathbf{C}_{p^n}$ is fully realizable if and only if $p^n = 1, 2, 3$ or $4$.
\end{proposition}
\begin{proof}
We have already proved that if $p$ is odd then $p^n = 1$ or $3$. For $p = 2$, we have already proved that $\mathbf{C}_2$ and $\mathbf{C}_4$ are fully realizable. For $3 \leq k \leq \infty$, $\mathbf{C}_{2^k}$ is not the group of units in a ring of characteristic 2 (\cite{CL15}).
\end{proof}

\begin{proposition}
If $T$ is a fully realizable torsion abelian group, then $T$ is a direct sum of finite cyclic groups.
\end{proposition}
\begin{proof}
Let $T$ be a fully realizable torsion abelian group. Recall that $T$ is the direct sum of its $p$-torsion subgroups $T_p$, over all primes $p$. Hence, each $T_p$ is fully realizable. Since an abelian group is divisible if and only if it is injective, the group $T_p$ is decomposable as the direct sum of a divisible group $D_p$ and a reduced group $R_p$ (reduced means the only divisible subgroup is 0). However, $D_p$ is a direct sum of copies of the quasi-cyclic group $\mathbf{C}_{p^\infty}$ (every divisible abelian group is a sum of copies of the rationals and quasicyclic groups, and $D_p$ is torsion), which is not realizable. So $D_p = 0$. If $R_p$ has unbounded exponent, then it has cyclic summands of arbitrarily high orders (\cite[{Exercise 1, p. 158}]{Fuchs15}), which is impossible by the previous proposition. So $R_p$ has bounded exponent and is therefore a direct sum of cyclic groups (\cite[Theorem 5.2]{Fuchs15}).
\end{proof}

In sum, a fully realizable torsion group $T$ must be a direct sum of fully realizable indecomposable cyclic groups ($\C_2, \C_3, \C_4$). The group $\C_3$ cannot appear more than once by Example \ref{c3xc3}. The group $\C_4$ cannot appear more than once by the following proposition (take $A = \C_4$).

\begin{proposition}\label{c4xc4} \label{Axc4}
If $A$ is a group that has $\mathbf{C}_4$ as a quotient, then $A \times \mathbf{C}_4$ is not fully realizable.
\end{proposition}
\begin{proof}
Let $A$ be a group such that there is a surjective group homomorphism $h\colon A \longrightarrow \C_4 = \langle y \rangle$. Assume to the contrary that there is a ring $R$ of characteristic 2, generated its units, that fully realizes $G = A \times \mathbf{C}_4.$ Pick any $x \in A$ such that $h(x) = y$. We have \[
\begin{aligned}
(1 + x + y)^2 &= 1 + x^2 + y^2 \\
(1 + x + y)^4 & = x^4
\end{aligned}
\] since $x$ and $y$ commute. Thus $1 + x + y$ is a unit and there is an element $z \in A$ and an integer $0 \leq b < 4$ such that
\begin{equation}
1 + x + y = z y^b.\tag{$\dag$}
\end{equation}
Equation ($\dag$) must hold after the application of any ring endomorphism induced by a group endomorphism. If we apply the group endomorphism $(s, t) \mapsto (s, 1)$ (fixing $x$ and $z$ and sending $y$ to 1), we obtain $x = z.$ If we apply the endomorphism $(s, t) \mapsto (1, t)$ (fixing $y$ and sending $x$ to 1), we obtain \[ 1 = y^{b-1}\] and hence $b = 1$. Now ($\dag$) is equivalent to the equation \[ (1 + x)(1+y) = 0.\] If we apply the endomorphism $(s, t) \mapsto (1, h(s)t)$ (sending $x$ to $y$ and fixing $y$; the map is a homomorphism since $h$ maps into an abelian group), we obtain \[ 0 = (1 + y)^2 = 1 + y^2,\] which implies $y$ has order 2. But $y$ has order 4, a contradiction.
\end{proof}

Write $G^I$ for $\underset{\scriptsize \alpha \in I}{\oplus} G$. Given that summands of fully realizable groups are fully realizable, to prove Theorem \ref{torsion} it remains to show that \[ \C_2^I \times \mathbf{C}_4 \times \mathbf{C}_3\] is fully realizable. Given Example \ref{2torsionxc3}, it suffices to prove that \[ \C_2^I \times \mathbf{C}_4\] is fully realizable. This follows from Proposition \ref{a24} below (take $V = 0$) (the torsion free summand $V$ is included in that proposition because we will later use it to classify the fully realizable finitely generated abelian groups).

Before proving Proposition \ref{a24}, we need to construct a ring that realizes (and, as we will later show, fully realizes) \[ \C_2^K = \underset{\alpha \in K}{\oplus} \langle x_\alpha \rangle,\] where $K$ is an arbitrary index set. For any subset $J = \{\alpha_1, \dots, \alpha_n\} \subseteq K$, write $x_J = x_{\alpha_1} \cdots x_{\alpha_n}.$ If $J$ is empty, $x_J = 1$.

\begin{proposition} \label{sumc2}
Let $K$ be an index set, let \[G = \C_2^K = \underset{\alpha \in K}{\oplus} \langle x_\alpha \rangle,\] and let $I\subseteq \mathbf{F}_2[G]$ denote the ideal generated by \[1 + x_\alpha + x_\beta + x_\alpha x_\beta\] for all $\alpha, \beta\in K$. Then, \[1 + x_A + x_B + x_A x_B \in I\] for all finite subsets $A, B \subseteq K$ and $\F_2[G]/I$ has unit group $G$.
\end{proposition}
\begin{proof} Let $R = \F_2[G]/I$, where $G$ and $I$ are defined in the statement of the proposition. All computations below take place in $R$. By induction on $|J|$, \[x_J = \sum_{\alpha \in J} x_\alpha + |J| + 1.\] Now:
\[
\begin{aligned}
x_A + x_B &= \sum_{\alpha \in A} x_\alpha + |A| + 1 + \sum_{\beta \in B} x_\beta + |B| + 1\\
&= \sum_{\gamma \in A\cup B \setminus (A \cap B)} x_\gamma + |A \cup B \setminus (A \cap B)|\\
&= x_{A\cup B \setminus (A \cap B)} + 1 \\
&= x_A x_B + 1.
\end{aligned}
\] Every element $z$ in $R$ may be uniquely written as \[z = a + \sum_{i=1}^n x_{\alpha_i}\] where the $\alpha_i$ are distinct and $a \in \{0, 1\}$. The element $z$ is a unit iff $a + n$ is odd (because $z^2 = a^2 + n$), iff $z = x_J$ where $J = \{\alpha_1, \dots \alpha_n\}$. This proves that $R^\times = G.$
\end{proof}

\begin{proposition} \label{a24}
Let $V$ be a torsion-free abelian group and let $K$ be an index set. 
\begin{enumerate}
    \item The group $G = V \times \C_2^K \times \mathbf{C}_4$ is fully realizable if and only if $V$ does not have $\C_4$ as a quotient. \label{yesc4}
    \item The group $V \times \C_2^K$ is fully realizable. \label{noc4}
\end{enumerate}
\end{proposition}
\begin{proof} We begin with item (\ref{yesc4}). By Proposition \ref{Axc4}, $G$ cannot be fully realizable if $V$ has $\C_4$ as a quotient. Hence, we need only prove the `if' direction of this proposition. Assume $V$ does not have $\C_4$ as a quotient.

Note that $\F_2[V]^\times = V$. If we write $\F_2[\C_4] = \F_2[y]/(y^4 + 1),$ then $(\F_2[\C_4]/(1 + y + y^2 + y^3))^\times = \C_4.$ Combining the work above in Propositions \ref{kgproduct} and \ref{sumc2}, we have that $(\F_2[G]/W)^\times = G$ where $W$ is the ideal generated by
\begin{enumerate}
\item[(1)] $vx_Jy^r + v + x_J + y^r$ for $v \in V, J \subseteq K, r \in \Z$;
\item[(2)] $1 + x_A + x_B + x_A x_B$ for finite subsets $A, B \subseteq K$; and
\item[(3)] $1 + y + y^2 + y^3$.
\end{enumerate}
Using (3), one can check that $1 + y^r + y^s + y^{r+s} \in W$ if and only if $r$ or $s$ is even (i.e., if at most one of $y^r$ and $y^s$ has order 4). By induction on $n$ and using the list above, it follows that in fact $W$ is the ideal generated by 
\begin{equation}
\prod_{i=1}^n u_i + \sum_{i=1}^nu_i + n + 1 \tag{$\star$}
\end{equation}
for all units $u_1, \dots, u_n \in G$ where at most one of the units $u_i$ has infinite order and at most one of the units $u_i$ has order 4 modulo $V$.

It remains to verify the condition on $W$ given in Proposition \ref{ideal-pres}; we must check that for any $\phi \in \mathrm{End}(G),$ if at most one of $u_1, \dots, u_n$ has infinite order and at most one unit has order 4 modulo $V$, then the same is true of the list of units $\phi(u_1), \dots, \phi(u_n)$. It is clear that at most one of $\phi(u_1), \dots, \phi(u_n)$ has infinite order. To see why at most one has order 4 modulo $V$, assume to the contrary that both $u_1 = v x_J y^a$ and $u_2 = w x_L y^b$ map to elements that have order 4 modulo $V$. Note that neither $\phi(v)$ nor $\phi(w)$ may have order 4 modulo $V$, for otherwise the inclusion of $V$ into $G$, followed by $\phi$, followed by projection onto $\C_4$ would be a surjective group homomorphism. Further, $\phi(x_J)$ and $\phi(x_L)$ have order at most 2 modulo $V$. This means that $\phi(y)^a$ and $\phi(y)^b$ must both have order 4 modulo $V$. This forces both $a$ and $b$ to be odd, and hence $u_1$ and $u_2$ both have order 4 modulo $V$, a contradiction.

One can delete the $\C_4$ factor in the argument above to prove item (\ref{noc4}). The difference is that the ideal is generated only by the generators from (1) (replace $y^r$ with 1) and (2) above, and in the characterization of the ideal as the ideal generated by the expressions ($\star$) we drop the requirement that at most one unit has order 4 modulo $V$ (all units will have order at most 2 modulo $V$).
\end{proof}

\section{Finitely generated abelian groups}

We will now prove the following theorem.

\begin{theorem} \label{fgag}
The fully realizable finitely generated abelian groups are the groups $W \times H$ and $W \times \C_\infty^n$ ($n \geq 1$) where $W$ is a finite elementary abelian $2$-group and $H$ is a subgroup of $\C_{12}.$
\end{theorem}

\noindent In light of Proposition \ref{a24} and Theorem \ref{torsion}, the theorem follows from the next proposition (take $A = \C_\infty$).

\begin{proposition}
If a group $A$ is a direct sum of indecomposable groups (e.g., if $A$ is finite, free abelian, or torsion-free abelian of finite rank) and has $\C_3$ as a quotient, then $A \times \C_3$ is not fully realizable. \label{c3iq}
\end{proposition}
\begin{proof}
Note that if $A$ has $\C_3$ as a quotient, then one of its indecomposable summands also has $\C_3$ as a quotient. Further, since summands of fully realizable groups are fully realizable, it suffices to prove the claim assuming $A$ is an indecomposable group with $\C_3$ as a quotient.

Suppose that $R,$ necessarily of characteristic 2, is generated by its units and fully realizes $R^\times = A \times \C_3$. Note that \[\begin{aligned}
\F_2[A \times \C_3] & \cong \F_2[A] \otimes_{\F_2} \F_2[\C_3] \\
& \cong \F_2[A] \otimes_{\F_2} (\F_2 \times \F_4) \\
& \cong \F_2[A] \times \F_4[A].
\end{aligned}\] Consider the map \[ \F_2[A] \times \F_4[A] \overset{\cong}{\longrightarrow} \F_2[A \times \C_3] \longrightarrow R,\] where the first map is determined by the chain of isomorphisms above and the second map is induced by the inclusion of $A \times \C_3$ into $R$. The first map maps $\{1\} \times \F_4^\times$ isomorphically onto $\{1\} \times \C_3 \subseteq R$, so $R$ has the form \[ \F_2[A]/I \times \F_4[A]/J\] where $J \neq \F_4[A]$. In the above ring direct product decomposition of $R$, let $S$ denote the first direct factor and let $T$ denote the second. 

Since $T^\times$ is a subgroup of $A \times \C_3$ that contains the direct factor $\C_3$, we have $T^\times = G \times \C_3$ for some group $G$. Thus, $(S\times T)^\times = S^\times \times \C_3 \times G$, and  $S^\times \times G = A$ since finite abelian groups have the cancellation property. But $A$ is indecomposable, so either $S^\times$ or $G$ is trivial.

If $G$ is trivial, then \[R = \F_2[A]/I \times \F_4\] (and the group of units of the first factor is $A$). Since there is a nontrivial map from $A$ to $\C_3$, this ring does not realize its group of units by Proposition \ref{products-fail}.

If $S^\times$ is trivial, then since $S^\times$ is generated by its units we have that either $S = \F_2$ or $S$ is the zero ring. In the latter case, we would have $R = \F_4[A]/J$; however, every ring endomorphism of the $\F_4$-algebra $\F_4[A]/J$ restricts to an injection on $\F_4$, making the quotient map $A \times \C_3 \longrightarrow A$ not realizable. So we may assume \[R = \F_2 \times \F_4[A]/J,\] where the group of units of the second factor is $A \times \C_3$. If $M$ is an endomorphism of $R$, then $M = (M_1, M_2)$ where \[M_1 \colon \F_2 \times \F_4[A]/J \longrightarrow \F_2\] and \[M_2 \colon \F_2 \times \F_4[A]/J \longrightarrow \F_4[A]/J.\] We must have $M_1(a,b) = a$, and the kernel of $M_2$ is either (1) $\F_2 \times K$ or (2) $\{0\} \times K$ where $K$ is an ideal of $\F_4[A]/J$. If $K = \F_4[A]/J$, then we're in case (2) and $M_2(a,b) = a$. Now suppose $K \neq \F_4[A]/J$.

In case (1), $M_2$ factors as \[ \xymatrix{M_2 \colon \F_2 \times \F_4[A]/J \ar@{>>}[r] & \F_4[A]/J\, \ar@{>>}[r]^-{\pi}& (\F_4[A]/J)/K \neq 0.}\] In this case $M_2(a,b) = \pi(b)$ and $M_2|_{\C_3} = \pi|_{\C_3}$ is injective because $K \neq \F_4[A]/J$. In case (2), $M_2$ embeds $\F_2 \times \F_4$ as a subring of the target, so again $M_2|_{\C_3}$ is injective.

In sum, $M|_{\C_3}$ is injective unless $M^\times = 1$. Thus, the map $(v, y) \mapsto (v, 1)$ is not realizable.
\end{proof}

The following theorem gives the consequent characterization of the fully realizable cyclic groups.

\begin{theorem}
\label{cyclicthm}
For $n \leq \infty$, $\mathbf{C}_n$ is fully realizable if and only if $n$ is a divisor of $12$ or $n = \infty$.
\end{theorem}

We include one last proposition concerning torsion-free abelian groups that do not have $\C_3$ (or, neither $\C_3$ nor $\C_4$) as a quotient (cf. Proposition \ref{c3iq}). It applies, for example, to $V = \Q$ (the rational numbers).

\begin{proposition} \label{vsg3} Let $V$ be a torsion-free abelian group and let $K$ be an index set.
\begin{enumerate}
    \item If $\C_3$ is not a quotient of $V$, then $V \times \C_2^K \times \C_3$ is fully realizable. \label{c3nq}
    \item If neither $\C_3$ nor $\C_4$ is a quotient of $V$, then $V \times \C_2^K \times \C_4 \times \C_3$ is fully realizable. \label{c3c4nq}
\end{enumerate}
\end{proposition}

\begin{proof} Let $V$ be a torsion-free abelian group and let $K$ be an index set. First, we prove item (\ref{c3nq}). If $\C_3$ is not a quotient of $V$, then there are no nontrivial maps from the fully realizable group $V \times \C_2^K$ to $\C_3$ or vice-versa. Thus, by Proposition \ref{products}, $V \times \C_2^K \times \C_3$ is fully realizable. The same argument may be used to establish item (\ref{c3c4nq}), because $V \times \C_2^K \times \C_4$ is fully realizable when $\C_4$ is not a quotient of $V$.
\end{proof}

\printbibliography
\end{document}